\documentclass[12pt]{article}
\usepackage[utf8]{inputenc} 
\usepackage{graphicx}
\usepackage{subfigure}
\usepackage{float}
\usepackage{amsmath}
\usepackage{amsfonts,dsfont,amssymb}

\usepackage{xcolor}

\usepackage{color}
\usepackage{ulem}

\newtheorem{theo}{Theorem}[section]
\newtheorem{defi}[theo]{Definition}
\newtheorem{remark}[theo]{Remark}

\newtheorem{prop}[theo]{Proposition}
\newcommand{\dem}{\textbf{Proof.}\ }
\begin{document}

\begin{center}
\vspace*{.7cm}
{\Large\bf A General Class of Fatigue-life Distributions of Birnbaum-Saunders Type}
\vspace*{1cm}
\\
 C.C.Y. Dorea$^{\footnotesize \mbox{a,}}$\footnote{ changdorea@gmail.com or
 changdorea@unb.br}
  R. Vila $^{\footnotesize \mbox{b,}}$\footnote{  rovig161@gmail.com or rovig161@unb.br } and
    F.S. Quintino$^{\footnotesize \mbox{a,}}$\footnote{Corresponding author: felipes.quintino2@gmail.com or  felipe.quintino@unb.br}
\\
 \vskip 6mm 
 $\phantom{.}^{\footnotesize \mbox{a}}$ 
Departament of Mathematics, Universidade de Brasilia,  70910-900 Brasilia, Brazil
\\
$\phantom{.}^{\footnotesize \mbox{b}}$ 
Departament of Statistics, Universidade de Brasilia,  70910-900 Brasilia, Brazil
\end{center}
\vskip 6mm

 \begin{abstract}
{\footnotesize
The derivation of Birnbaum-Saunders (BS) fatigue-life distribution is based on the asymptotic normality of crack damages' partial sums. We address the situation when this fails and the crack damages possess heavy-tailed distribution. In this note, a general class of fatigue-life distributions of BS type is proposed.
}
 \end{abstract}

	\section{Introduction}
Let $X_1, X_2, \cdots$ be a sequence of random variables (r.v.'s) each representing the cumulative fatigue damage over a cycle. It is assumed that $X_j's$ are non-negative, independent and identically distributed (i.i.d) with a common distribution $F_X$. The Birnbaum-Saunders (BS) fatigue-life model (cf. \cite{birnbaum1969new} and \cite{leiva2015birnbaum}) proposes a distribution $\rho_n$ the smallest number of cycles $N_*$ such that the partial sum $S_n = \sum_{j=1}^n X_j$ exceeds a given threshold value $s_*$.
The derivation of BS fatigue-life distribution is based on the asymptotic normal approximation $S_n\cong N(n\mu_X, n\sigma_X^2)$ where $\mu_X\neq0$ and $\sigma_X>0$ are, respectively, the mean and the variance of $F_X$. More specifically, let
\begin{equation*}
    N_* = \inf\left\{ n; S_n = \sum_{j=1}^n X_j >s_* \right\}
\end{equation*}
and let $T_*$ be the continuous counterpart of $N_*$. Then the distribution of $T_*$ satisfies,
\begin{equation}\label{eq_1.1}
    \mathbb{P}(T_*\leqslant t) 
    = 
    \Phi\left(\frac{1}{a}\, \xi\left(\frac{t}{b}\right) \right), ~~t>0,
\end{equation}
where $\Phi$ is the cumulative distribution function (c.d.f.) of $N(0,1)$,
\begin{equation}\label{eq_1.2}
    a=\frac{\sigma_X}{\sqrt{\mu_X s_*}},
    \quad  
    b=\frac{s_*}{\mu_X}
    \quad 
    \mbox{and}
    \quad 
    \xi(x) = \sqrt{x} - \frac{1}{\sqrt{x}}, \quad 
    x>0.
\end{equation}

Now, we assume that the crack damages $X_1, X_2, \cdots$ possess a heavy-tailed distribution in the sense that $\mathbb{E}(X_k^2)=\infty$. 
Then the asymptotic behavior of normalized partial sums cannot be Gaussian. In this case, the class of stable distributions $\{S_\alpha(\sigma, \beta, \mu); 0<\alpha\leqslant 2\}$ plays a central role, a similar role the normal distribution ($\alpha=2$) plays among distributions with finite second moments. 
A useful tool to handle stable laws is provided by the Mallows (or Wasserstein) distance $d_r(F,G)$ that measures the discrepancy between two distribution functions $F$ and $G$. 
For $r>0$ define
\begin{equation*}
    d_r(F,G) 
    = 
    \inf_{(X,Y)}\left\{\mathbb{E} |X-Y|^r\right\}^{1/r},  \quad X\overset{d}{=}F, ~Y\overset{d}{=}G,
\end{equation*}
where the infimum is taken over all random vectors $(X,Y)$ with marginal distributions $F$ and $G$, respectively, and $\overset{d}{=}$ denotes equality in distribution. Besides applications to several statistical methods such as bootstraps and goodness-of-fit tests (cf. \cite{bickel1981some} and \cite{del2005asymptotics}), this metric has also been successfully used to establish central limit theorems (CLT) type results for heavy-tailed distributions (see \cite{johnson2005central} or \cite{barbosa2009note}).

Our Proposition \ref{prop_2.3} shows that, for $1<\alpha<2$, if there exists an $\alpha$-stable distribution $G_\alpha$ such that $d_\alpha(F_X, G_\alpha)<\infty$, then
\begin{equation}\label{eq_1.4}
    S_n \cong S_\alpha\left(\sigma n^{1/\alpha}, 0, n\mu_X\right),
\end{equation}
where $\sigma>0$ is a constant. And this leads us to a general class of BS type fatigue-life distributions,
\begin{equation}\label{eq_1.5}
    \mathbb{P}(T_* \leqslant t) = \Phi_\alpha\left(\frac{1}{a}\, \xi_\alpha\left(\frac{t}{b}\right)\right), \quad t>0,
\end{equation}
where $\Phi_\alpha$ is the c.d.f. of $G_\alpha$ (cf. Proposition \ref{prop_2.4}).
Also, if $\alpha=2$ then (\ref{eq_1.5}) reproduces (\ref{eq_1.1}) and (\ref{eq_1.2}) (see Remark \ref{remark_2.5}).
For $G_\alpha=S_\alpha(\sigma, 0, \mu)$, estimation issues concerning $\alpha$ and $\sigma$ are also addressed, see, \eqref{eq_2.5} and \eqref{eq_2.6}.

\section{General Birnbaum-Saunders Distribution}
First, we state some properties of stable distributions. For an extensive treatment on the matter, we refer the reader to \cite{SamorodnitskyTaqqw94}.

\begin{defi}
    For $0<\alpha\leqslant 2$, we say that $S_\alpha(\sigma, \beta, \mu)$ is an $\alpha$-
    {\bf stable distribution} with scale parameter $\sigma>0$, skewness parameter $|\beta|\leqslant1$ and shift parameter $\mu\in\mathbb{R}$, if for any $n\geqslant2$, there are real numbers $d_n$ such that
    \begin{equation*}
    Y_1+\cdots+Y_n \overset{d}{=} n^{1/\alpha} Y+d_n,
    \quad  Y\overset{d}{=}S_\alpha(\sigma, \beta, \mu),
    \end{equation*}
    where $Y_1, Y_2, \cdots$ are independent copies of $Y$. If $\beta=0$ we say that $Y$ has a strictly $\alpha$-stable distribution.
\end{defi}

\begin{prop}\label{prop_2.2}
    Let $Y\overset{d}{=}S_\alpha(\sigma, \beta, \mu)$. Then
    \begin{itemize}
    \item[(a)] In general $\sigma>0$ does not represent the variance of the distribution. Except for Gaussian case, $S_2(\sigma, 0, \mu)=N(\mu, \sigma^2)$.
    
        \item[(b)] If $\alpha>1$ then $\mathbb{E}(Y)=\mu$ and $d_n=\mu(n-n^{1/\alpha})$.

        \item[(c)] If $0<\alpha'<\alpha<2$ then $\mathbb{E}|Y|^{\alpha'}<\infty$ and $\mathbb{E}|Y|^{\alpha}<\infty$.

        \item[(d)] If $\beta=0$ then for constant $a$ and $b$ we have $aY+b\overset{d}{=}S_\alpha(|a|\sigma, 0, a\mu+b)$.

        \item[(e)] The characteristic function of $S_\alpha(\sigma, 0, 0)$ is given by $\varphi_{S_\alpha}(t) = e^{-\sigma^2 t^2}$.
    \end{itemize}
\end{prop}

As for the Mallows distance, it is worth pointing out that, from the practical point of view, it is fairly simple to compute. The representation result from \cite{dorea2012conditions} shows that: for $r\geqslant1$, we have
\begin{equation}\label{eq_2.1}
    d_r^r(F, G) 
    = 
    \mathbb{E}|X^*-Y^*|^r = \int |x-y|^r {\rm d}(F(x)\wedge G(y)),
\end{equation}
where $X^*\overset{d}{=}F$, $Y^*\overset{d}{=}G$ and $(X^*, Y^*) \overset{d}{=}F\wedge G$. That is, 
\begin{equation*}
\mathbb{P}(X^*\leqslant x, Y^*\leqslant y) = \min\{F(x), G(x) \},
 \quad \forall (x,y)\in \mathbb{R}^2.   
\end{equation*}

\begin{prop}\label{prop_2.3}
    Let $1<\alpha<2$ and assume that $d_\alpha(F_X, G_\alpha)<\infty$, where $G_\alpha=S_\alpha(\sigma, 0, \mu)$. Then
    \begin{equation}\label{eq_2.2}
        \frac{S_n - n\mu_X}{n^{1/\alpha}} \overset{d}{\longrightarrow} S_\alpha(\sigma, 0, 0),
    \end{equation}
    and
    \begin{equation}\label{eq_2.3}
        S_n \overset{d}{=} S_\alpha(\sigma n^{1/\alpha}, 0, n\mu_X),
    \end{equation}
    with $\overset{d}{\longrightarrow}$ denoting convergence in distribution.
\end{prop}
\dem (a) Theorem 1 from \cite{barbosa2009note} shows that for a sequence of independent r.v.'s satisfying Lindeberg type conditions, we have
\begin{equation*}
    d_\alpha(F_n, G_\alpha) \underset{n}{\longrightarrow}0, \quad  
    F_n\overset{d}{=} \frac{S_n-n\mu_X + n^{1/\alpha}\mu}{n^{1/\alpha}}.
\end{equation*}
Since we are dealing with i.i.d. r.v.'s the conditions are trivially satisfied. Let $X\overset{d}{=}G_\alpha$ and take the joint distribution
\begin{equation*}
    \mathbb{P}\left( \frac{S_n-n\mu_X + n^{1/\alpha}\mu}{n^{1/\alpha}} \leqslant x, Z\leqslant z \right) = \min\left\{ F_n(x), G_\alpha(z) \right\}.
\end{equation*}
Using the representation (\ref{eq_2.1}) we can write
\begin{equation*}
    d_\alpha^\alpha(F_n, G_\alpha) 
    = 
    \mathbb{E}\left| \frac{S_n-n\mu_X + n^{1/\alpha}\mu}{n^{1/\alpha}} - Z \right|^\alpha   \underset{n}{\longrightarrow} 0.
\end{equation*}
From Proposition \ref{prop_2.2}, we have
\begin{equation*}
    \mathbb{E}\left| \frac{S_n-n\mu_X}{n^{1/\alpha}} - Z' \right|^\alpha \underset{n}{\longrightarrow} 0,
\end{equation*}
where $Z'=Z-\mu \overset{d}{=}S_\alpha(\sigma, 0, 0)$. And (\ref{eq_2.2}) follows from the moment convergence of order $\alpha$.

(b) Since $({S_n - n\mu_X})/({\sigma n^{1/\alpha}})\overset{d}{=}S_\alpha(1,0,0)$, we have
\begin{eqnarray}
    \nonumber \mathbb{P}(N_*\leqslant n) &=& \mathbb{P}(S_n > s_*)\\[0,2cm]
    \nonumber &=& \mathbb{P}\left(\frac{1}{\sigma}\frac{S_n - n\mu_X}{n^{1/\alpha}} > \frac{1}{\sigma}\frac{s_* - n\mu_X}{n^{1/\alpha}} \right)\\[0,2cm]
    \nonumber &\cong& 1-\Phi_\alpha\left(\frac{1}{\sigma}\frac{s_* - n\mu_X}{n^{1/\alpha}}\right),
\end{eqnarray}
where $\Phi_\alpha$ is the c.d.f. of $S_\alpha(1.0,0)$. And \eqref{eq_2.3} follows readily.
\begin{flushright}
    $\square$
\end{flushright}

\begin{prop}\label{prop_2.4}
    Let $1<\alpha\leqslant 2$ and assume $d_\alpha(F_X, G_\alpha)<\infty$, where $G_\alpha=S_\alpha(\sigma, 0, \mu)$. Then a general class of fatigue-life distributions of Birnbaum-Saunders type can be generated by the relation
    \begin{equation}\label{eq_2.4}
        \mathbb{P}(T_*\leqslant t) = \Phi_\alpha\left( \frac{1}{a_\alpha}\, \xi_\alpha\left(\frac{t}{b_\alpha}\right) \right), \quad t>0,
    \end{equation}
    where $\Phi_\alpha$ is the c.d.f. of $S_\alpha(\sigma, 0, 0)$, 
    \begin{align*}
    a_\alpha = \frac{1}{\mu_X^{1/\alpha} s_*^{1-1/\alpha}}, 
    \quad 
    b_\alpha = \frac{s_*}{\mu_X} 
    \quad
    \text{and} 
    \quad 
    \xi_\alpha(x) = x^{1-1/\alpha} - \frac{1}{x^{1/\alpha}}, 
    \quad x>0. 
    \end{align*}
\end{prop}
\dem (a) For $1<\alpha<2$, since (\ref{eq_2.2}) holds, we have
\begin{eqnarray}
    \nonumber \mathbb{P}(N_*\leqslant n) &=& \mathbb{P}\left( \frac{S_n - n\mu_X}{n^{1/\alpha}}> \frac{s_* - n\mu_X}{n^{1/\alpha}} \right)\\[0,2cm]
    \nonumber &\cong& 1- \Phi_\alpha\left(\frac{s_* - n\mu_X}{n^{1/\alpha}}\right)\\[0,2cm]
    \nonumber &=& \Phi_\alpha\left( \frac{n\mu_X - s_*}{n^{1/\alpha}} \right).
\end{eqnarray}
Set the reparametrization $a_\alpha = {1}/({\mu_X^{1/\alpha} s_*^{1-1/\alpha}})$ and $b_\alpha = {s_*}/{\mu_X}$. And write
\begin{eqnarray}
    \nonumber \frac{n\mu_X - s_*}{n^{1/\alpha}} &=& \mu_X^{1/\alpha} s_*^{1-1/\alpha} \left[ n^{1-1/\alpha} \left(\frac{\mu_X}{s_*}\right)^{1-1/\alpha} - \frac{1}{n^{1/\alpha}}\left(\frac{s_*}{\mu_X}\right)^{1/\alpha}\right]\\[0,2cm]
    \nonumber &=& \frac{1}{a_\alpha} \left[ n^{1-1/\alpha} \left(\frac{1}{b_\alpha}\right)^{1-1/\alpha} - \frac{1}{n^{1/\alpha}}b_\alpha^{1/\alpha}\right]\\[0,2cm]
     \nonumber &=& \frac{1}{a_\alpha} \left[  \left(\frac{n}{b_\alpha}\right)^{1-1/\alpha} - \left(\frac{b_\alpha}{n}\right)^{1/\alpha}\right]\\[0,2cm]
      \nonumber &=& \frac{1}{a_\alpha} \xi_\alpha\left(\frac{n}{b_\alpha}\right).
\end{eqnarray}
For (\ref{eq_2.4}), let $T_*$ and $t$ be the continuous counterpart of $N_*$ and $n$.

(b) Let $\alpha=2$. Since $S_2(\sigma, 0, \mu) = N(\mu, 2\sigma^2)$, $G_2$ possesses moments of all order. Condition $d_2(F_X, G_2)<\infty$ assures that $F_X$ has finite variance. Then
\begin{equation*}
    \frac{S_n-n\mu_X}{\sqrt{n}} \overset{d}{\longrightarrow} N(0, \sigma_X^2) = S_2\left(\frac{\sigma_X}{\sqrt{2}}, 0, 0\right).
\end{equation*}
In this case, we necessarily have $\Phi_2$ the c.d.f. of $S_2(\sigma, 0, 0)$ where $\sigma={\sigma_X}/{\sqrt{2}}$.
Exactly the same proof as in (a) shows (\ref{eq_2.4}) holds for $\alpha=2$.
\begin{flushright}
    $\square$
\end{flushright}

\begin{remark}\label{remark_2.5}
    \begin{itemize}
        \item[(a)] It is worth pointing out that \eqref{eq_2.4} for $\alpha=2$ coincides with \eqref{eq_1.1} and (\ref{eq_1.2}). Setting $\alpha=2$ we get $\xi_2=\xi$, $a_2={1}/{\sqrt{\mu_X s_*}}$ and $b_2={s_*}/{\mu_X}$. Now,
        \begin{equation*}
            \Phi_2(y) =\Phi_{S_2\left(\frac{\sigma_X}{\sqrt{2}}, 0, 0\right)}(y) = \Phi_{N(0,\sigma_X^2)}(y) = \Phi_{N(0,1)}(\sigma_X^{-1} y)
        \end{equation*}
        and 
        \begin{equation*}
            \Phi_2\left(\frac{1}{a_2}\, \xi_2\left(\frac{t}{b_2}\right)\right) = \Phi_{N(0,1)}\left(\frac{1}{\sigma_X a_2}\, \xi_2\left(\frac{t}{b_2}\right) \right).
        \end{equation*}
        From \eqref{eq_1.2}, we have $a=\sigma_X a_2$ and $b=b_2$.

        \item[(b)] Unlike the Gaussian case, unfortunately for $1<\alpha<2$ there is no close form for the c.d.f. $\Phi_\alpha$.
        However, computationally it is interesting to express (\ref{eq_2.4}) in terms of $H$-functions (see \cite{rathie2017exact} and the references therein). 
It follows from reference \cite[Proposition 1]{rathie2017exact} that 
\begin{equation*}
    \mathbb{P}(T_*\leqslant t) = \frac{1}{2} +  \frac{1}{a_\alpha \sigma} \xi_\alpha\left(\frac{t}{b_\alpha}\right) H_{1,2}^{3,3}\left[ \left|\frac{1}{a_\alpha \sigma} \xi_\alpha\left(\frac{t}{b_\alpha}\right)\right|\left| \begin{array}{c}
     (1-1/\alpha, 1/\alpha), (0,1), (1/2, 1/2)   \\
       (0, 1),  (1/2, 1/2), (-1,1)
  \end{array}  \right.\right],
\end{equation*}
where the $H$-function (cf. \cite{MathaiSaxenaHaubold10}) is defined by
\begin{equation*}
     H_{p,q}^{m,n}\left[ z\Big| \begin{array}{c}
     (a_1, A_1), \cdots, (a_p, A_p)   \\
       (b_1, B_1), \cdots, (b_q, B_q) 
  \end{array}  \right] = \frac{1}{2\pi i} \int_L   \frac{ \prod_{k=1}^m \Gamma (b_j + B_js) \prod_{j=1}^n \Gamma (1-a_j - A_js)}{\prod_{k=m+1}^q \Gamma (1-b_j - B_js) \prod_{j=n+1}^p \Gamma (a_j + A_js) } z^{-s} {\rm d}s,
\end{equation*}
in which  $0 \leqslant m \leqslant q$, $0 \leqslant n \leqslant p$ (not both $m$ and $n$ simultaneously zero), $A_j > 0$, $j=1,\cdots,p$,  $B_k > 0$, $k=1,\cdots,q$, $a_j$ and $b_k$ are complex numbers such that no poles of $\Gamma (b_k + B_k s)$, $k=1,\cdots,m,$ coincide with poles of $\Gamma (1-a_j - A_j s)$, $j=1,\cdots, n$, with $\Gamma$ being the complete gamma function. 
$L$ is a suitable contour $w -i \infty$ to $w + i \infty$, 
$w\in\mathbb{R}$, separating the poles of the two types mentioned above. For more details, see \cite{MathaiSaxenaHaubold10}.

 
    \end{itemize}
\end{remark}

For the approximation (\ref{eq_1.4}) the stability and scaling parameters, $\alpha$ and $\sigma$, remain to be estimated. Using the well-known Hill's estimate (cf. \cite{hill1975simple}), we have
\begin{equation}\label{eq_2.5}
    \frac{1}{k_n} \sum_{j=n-k_n+1}^n \log \frac{X_{(j)}'}{X_{(n-k_n)}'} \overset{p}{\longrightarrow} \frac{1}{\alpha},
\end{equation}
where $\overset{p}{\longrightarrow}$ denotes convergence in probability, $k_n\uparrow\infty$, ${k_n}/{n}\underset{n}{\rightarrow}0$, $X_j' = X_j - \mu_X$ and 
\begin{equation*}
    X_{(1)}' \leqslant X_{(2)}' \leqslant \cdots \leqslant X_{(n)}'
\end{equation*}
 are the order statistics of $(X_1', X_2', \cdots, X_n')$.

As for the scaling parameter $\sigma$, we make use of the results from \cite[Corollary 1]{dorea2007levy}. Consider the characteristic function
\begin{equation*}
    \varphi_Y(t) = \mathbb{E}\{\exp({itY})\} = \exp({-\sigma^\alpha |t|^\alpha}), \quad  Y\overset{d}{=} S_\alpha(\sigma, 0, 0).
\end{equation*}
The inverse transform of $\varphi_Y$ is given by
\begin{equation*}
    l_{\alpha,\sigma}(x) = \frac{1}{\pi} \int_0^\infty \cos(xt) \exp({-\sigma^\alpha |t|^\alpha}) {\rm d}t
\end{equation*}
and its value at point $x=0$,
\begin{equation*}
    l_{\alpha, \sigma}(0) 
    = \frac{1}{\pi \alpha\sigma}\,  \Gamma\left(\frac{1}{\alpha}\right).
\end{equation*}

Let $\varepsilon_n\downarrow0$ and define
\begin{equation*}
    Y_j^{(n)} = \frac{X_{(j-1)k_n +1}'+\cdots+X_{j k_n}'}{k_n^{1/\alpha}}
\end{equation*}
and
\begin{equation*}
    \hat{l}_n(0) = \frac{1}{2\varepsilon_n} \frac{1}{r_n} \sum_{j=1}^{r_n} \mathds{1}(|Y_j^{(n)}|\leqslant \varepsilon_n),
\end{equation*}
where $r_n=[{n}/{k_n}]$, the largest integer and $\mathds{1}(\cdot)$, the indicator function, then besides \eqref{eq_2.5} we also have 
\begin{equation}\label{eq_2.6}
\hat{\sigma}_n = \frac{\Gamma(1/\alpha)}{\pi \alpha} 
\big[\hat{l}_n(0)\big]^{-1}\overset{p}{\longrightarrow}\sigma.
\end{equation}

  \bibliographystyle{abbrv}

\end{document}